\documentclass[10pt,reqno]{amsart}
\usepackage{amsmath, amsthm}

\makeatletter

\numberwithin{equation}{section}
\newtheorem{thm}{\indent\bf Theorem}[section]
\newtheorem{lemma} [thm] {\indent\bf Lemma}

\newtheorem{Rem}   [thm] {\indent\bf Remark}

\begin{document}

\title[Embedding result]{An embedding result}
\author{A. Canale}
\address{Dipartimento di Matematica, Universit\`a degli Studi di Salerno, 
Via Giovanni Paolo II n. 132, 84084 FISCIANO (Sa), Italy.}
\email{acanale@unisa.it}

\begin{abstract}

In unbounded subset $\Omega$ in $R^n$ we study the operator $u\rightarrow gu$ 
as an operator defined in the Sobolev space $W^{r,p}(\Omega)$
and which takes values in $L^p(\Omega)$. The functions $g$ belong to wider 
spaces of $L^p$ connected with the Morrey type spaces.
The main result is an embedding theorem from which we 
can deduce a Fefferman type inequality.

\end{abstract}

\maketitle
\thispagestyle{empty}

\bigskip

{\bf AMS Subject Classifications}:  35J25, 46E35

\bigskip

{\bf Key Words}: elliptic equations, multiplication operator, embedding theorems.

\bigskip\bigskip\bigskip

\section {Introduction} 

\bigskip 
Let $\Omega$ be an unbounded open subset in $R^n$. 

In literature there are different results about the study of 
{\it multiplication} operator for a suitable function 
$g:\;\Omega\rightarrow C$ 
\begin{equation}\label{multiplication}
u\longrightarrow gu,\
\end{equation}
as an operator defined in a Sobolev space (with or without weight) 
and which takes values in a $L^p(\Omega)$ space.

In $W^{1,p}_0 (\Omega)$ or in $W^{1,p}(\Omega)$ 
with $\Omega$ regular enough, 
reference results are 
some well-known inequalities which state the boundedness of (\ref{multiplication}): 
Hardy type inequalities (see 
H.Brezis \cite{2}, A.Kufner \cite{11}, J.Ne{\v c}as \cite{12}) 
when $g(x)$ is an appropriate power of the distance of $x$ from a subset 
of $\partial \Omega$, C.Fefferman inequality \cite{10} (see, e.g. 
F.Chiarenza-M.Franciosi \cite{8}, F.Chiarenza-M.Frasca \cite{9}) 
obtained when $g$ belongs to a suitable Morrey space. 

In this paper we study the operator (\ref{multiplication}) 
in the Sobolev space $W^{r,p}(\Omega)$, 
$r\in N$,$1\le p\le+\infty$ when $g$ belongs to suitable spaces
$S^{p,s}$ defined in Section 3.

One of the aspects of our interest in this type of inequality lies in the fact that
the embedding results are useful tools to prove a priori bounds
when studying elliptic equations. For applications in the study of the a priori bounds
see \cite{3}, \cite{4}, \cite{5}, \cite{6}, \cite{7}. The spaces considered in
some of these papers are connected to the spaces $S^{p,s}$.

These spaces are wider than $L^p$ spaces, than classical Morrey space and 
are connected, for suitable values of $s$,
to the well known Morrey spaces defined 
when $\Omega$ is an unbounded open subset in $R^n$.
   
We characterize the classes of functions in $S^{p,s}$ and their
inclusion properties in Section 3.
In Section 4 we state the main result of this paper.
We prove that if $g\in S^{q,s}(\Omega)$ for an appropriate 
$q\in [p,+\infty[$ and $s\ge 0$, then (\ref{multiplication})
defines a bounded operator from $W^{r,p}(\Omega)$ 
in $L^p(\Omega)$. 

We remark that our results imply that, for any function 
$g:\;\Omega \rightarrow C$ such that 
$$
\sup_{x\in \Omega\atop
\rho\in]0,d]}\rho^{q-n}\int_{\Omega\cap 
B(x,\rho)}|g|^q<+\infty,\quad d>0,
\quad 
q\ge p,\;q\ge n,\;q>n\quad\hbox{se}\quad n=p>1,
$$ 
the following inequality holds
$$
\|gu\|_{L^p(\Omega)}\le c 
\|g\|_{S^{q,\frac{s}{p}}(\Omega)} 
\|\nabla u\|_{L^p(\Omega)} 
\qquad \forall u\in W^{1,p}_0(\Omega).
$$

\bigskip\bigskip\bigskip

\section {Notations}

\bigskip 

Let $R^n$ be the $n$-dimensional real euclidean space. We set 
$$B(x,r)=\{y\in R^n :\; |y-x|< r\},\quad B_r=B(0,r)\, 
\qquad \forall x\in R^n,\>\forall r\in R_+.$$ 

For any $x\in R^n$, we call {\it 
open infinite cone} having vertex at $x$ every set of the type 
$$
\{x+\lambda (y-x):\,\lambda\in 
R_+,\,\,|y-z|<r\},
$$ 
where $r\in R_+$ and $z\in R^n$ are such that $|z-x|>r$. 

For all $\theta\in ]0,\pi/2[$ and for all $x\in R^n$ 
we denote by $C_\theta(x)$ an open infinite cone 
having vertex at $x$ and opening $\theta$. 

For a fixed $C_\theta(x)$, we set 
$$
C_\theta(x,h)=C_\theta(x)\cap B(x,h)\,,\qquad\forall h\in R_+.
$$ 
Let $\Omega$ be an open set in $R^n$. 
We denote by $\Gamma (\Omega,\theta,h)$ the family of open cones 
$C\subset\subset\Omega$ of opening $\theta$ and height $h$. 

We assume that the following hypothesis is satisfied: 

\begin{itemize}
\medskip

\item[$h_1$)]  There exists 
$\theta\in ]0,{\pi/ 2}[$ such that 
$$\forall x\in \Omega\qquad \exists C_{\theta}(x)\quad \hbox{such that} \quad 
\overline 
{C_{\theta}(x, \rho)}\subset \Omega.
$$ 
\end{itemize}

\bigskip\bigskip 

\section { Spaces $S^{p,s}(\Omega)$} 

\bigskip 

Let $(\Omega(x))_{x\in \Omega}$ be the family of open sets 
in $R^n$ defined as
$$\Omega(x)= B(x,\rho)\cap \Omega, 
\qquad x\in \Omega,
\qquad \rho>0.$$ 
If $1\le p< +\infty$ and $s\in R$, we denote by $S^{p,s}(\Omega)$ 
the space of functions $g\in L^p_{loc}(\overline\Omega)$ such that 
\begin{equation}\label{space}
\|g\|_{S^{p,s}(\Omega)}=\sup_{x\in \Omega\atop \rho\in ]0,d]}\left(\rho^{s-n/p}\; 
\|g\|_{L^p(\Omega(x))}\right)<+\infty, 
\qquad d>0,
\end{equation} 
with the norm defined by (\ref{space}). 

We remark that 
$$
L^\infty(\Omega)\hookrightarrow S^{p,s}(\Omega) 
\qquad\forall p\in[1,+\infty[ 
\quad\hbox{and}\quad\forall s\ge 0.
$$ 
Indeed, if $g\in L^\infty(\Omega)$, we get
 
\begin{equation}
\begin{split}
\|g\|_{S^{p,s}(\Omega)}&= 
\sup_{x\in \Omega\atop \rho\in]0,d]}\left(\rho^{s-n/p}
\;\| g\|_{L^p(\Omega(x))}\right)\le
\\&
\le \| g\|_{L^{\infty}(\Omega)} 
\sup_{x\in\Omega\atop \rho+\in]0,d]}\rho^s 
\left(\rho^{-n/p} |\Omega(x)|^{1/p}\right) 
=c\|g\|_{L^\infty(\Omega)}
\end{split}
\end{equation}
where $c=c(n,p)$.

The following inclusions hold 
\begin{equation}\label{inclusions 1}
L^r(\Omega) 
\hookrightarrow S^{q,s}(\Omega), 
\qquad  s\ge \frac{n}{q},
\qquad 1\le q\le r\le +\infty,
\end{equation}
and
\begin{equation}\label{inclusions 2}
S^{q,s}(\Omega) 
\hookrightarrow 
S^{p,s}(\Omega)\qquad \quad 1\le p\le q\le +\infty,
\end{equation}

\noindent In particular
\begin{equation}\label{inclusion}
L^p(\Omega) 
\hookrightarrow S^{p,\frac{n}{p}}(\Omega).
\end{equation}

\noindent Indeed 
\begin{equation}
\begin{split}
\| g\|_{S^{q,s}(\Omega)} &\le
\sup_{x\in \Omega\atop \rho\in]0,d]}(\rho^{s-{n\over q}}
\|g\|_{L^r(\Omega(x))}|\Omega(x)|^{{1\over q}-{1\over r}})\le
\\&
\le c_1\sup_{x\in \Omega\atop \rho\in]0,d]}
\rho^{s-{n\over r}}\|g\|_{L^r(\Omega(x))}
\le c_1
\|g\|_{L^r(\Omega)}, 
\end{split}
\end{equation}
where $c_1$ is a constant independent of $g$. 
From which we deduce also that
$$
\| g\|_{S^{p,s}(\Omega)} \leq c_2 
\sup_{x\in \Omega\atop \rho\in]0,d]}
\rho^{s-{n\over q}} \|g\|_{L^q(\Omega(x))},
$$ 
where $c_2$ is a constant independent of $g$.

Regard to the inclusion (\ref{inclusion}) we note that, for example, 
the constant functions belong to $S^{p,\frac{n}{p}}(\Omega)$ and
do not belong to $L^p(\Omega)$. Furthermore the function 
$\frac{1}{1+|x|^\alpha}$ belongs to $S^{p,\frac{n}{p}}$ for any $\alpha>0$
but does not belong to $L^p$ if $\alpha\in ]0,\frac{n}{p}[$. 

\bigskip

\begin{Rem} When $\Omega=R^n$ the space $S^{p,s}$ includes 
the classical Morrey spaces $L^{p,n-sp}$, $0\le s\le \frac{n}{p}$, defined as
the space of functions $g\in L^p_{loc}(R^n)$ such that 
\begin{equation}
\|g\|_{L^{p,n-sp}(R^n)}=\sup_{x\in \Omega\atop \rho>0}
\left(\rho^{s-n/p}\; \|g\|_{L^p(B(x,\rho)}\right)<+\infty.
\end{equation}
\end{Rem}

\bigskip\bigskip\bigskip

\section{Embedding result}

\bigskip

Let us consider the function 
\begin{equation}\label{phi}
\phi: \;(x,y)\in \Omega\times \Omega 
\longrightarrow
\begin{cases} 
1 & \hbox{if}\quad y\in \Omega(x)\\
0&\hbox{if}\quad y\not\in \Omega(x).
\end{cases}
\end{equation}
and, for any $x\in \Omega$, we set 
\begin{equation}\label{E}
E(x)=\{y\in \Omega : x\in \Omega(y)\}.
\end{equation}

\bigskip\bigskip

\begin{lemma} For any $x\in \Omega$, $E(x)$ is a 
measurable set and there exist $c_1,c_2\in R_+$ such that 
\begin{equation}\label{mis E}
c_1\rho^n\leq |E(x)|\le c_2\rho^n\qquad \forall x\in 
\Omega.
\end{equation}
\end{lemma}

\medskip 

\begin{proof} Clearly the function $\phi$ defined by (\ref{phi}) is a 
measurable function. Then, for any fixed $ y\in \Omega$, 
the function $\phi^y: x\in \Omega \rightarrow \phi(x,y)$ 
is measurable. Since $\phi^y$ is the characteristic function of $E(y)$, we have 
that $E(y)$ is measurable. 

Now we prove that (\ref{mis E}) holds. 
The inequality on the right is easily proved.
We will prove the inequality on the left. 

Let us consider 
$C_{\theta}(x,\rho)$ such that 
$\overline {C_{\theta}(x,\rho)}\subset \Omega$. 
We get
$$
C_{\theta}(x,\rho)\subset E(x).
$$ 
In fact let $y\in C_{\theta}(x,\rho)$. Then there exists a cone 
$C\in\Gamma (\Omega,\theta,h)$ such that $x,y\in C$.  
So
$$
x\in B(y,\rho)\cap \Omega \Longrightarrow
y\in E(x).
$$
Thus the inequality (\ref{mis E}) is stated.

\end{proof}

\bigskip
 
Now we state a Lemma which we will use in the proof of the embedding result.

\bigskip\bigskip

\begin{lemma} If $h_1)$ holds, then, 
for any $s\ge 0$, we have $v\in L^1(\Omega)$ if and only if the map 
$x\in \Omega\rightarrow \rho^{-n}|v|_{L^1(\Omega(x))}$ belongs to 
$L^1(\Omega)$. 
Therefore there exist $c_1,c_2\in R_+$ such that 
\begin{equation}\label{map}
c_1\|v\|_{L^1(\Omega)}\le 
\int_{\Omega} \rho^{-n}\|v\|_{L^1(\Omega(x))}dx\le c_2\|v\|_{L^1(\Omega)}\qquad 
\forall v\in L^1(\Omega).
\end{equation}
\end{lemma}

\medskip 

\begin{proof}
The result is a consequence of the relation 
\begin{equation}
\begin{split}
\int_{\Omega}\rho^{-n}\|v\|_{L^1(\Omega(x))}dx&= 
\int_{\Omega} \rho^{-n}\int_{\Omega}|v(y)| \phi(x,y)dxdy= 
\\&
=\int_{\Omega}|v(y)|dy \int_{E(y)}\rho^{-n}dx 
\end{split}
\end{equation}
and of the Lemma 4.1.
\end{proof}

\bigskip
 
Let $r, s, p, q$ be real number with the condition 

$$
h_2)\quad r\in N,\quad s\ge 0, \qquad 
1\le p\le q<+\infty, \qquad q\ge 
{n\over r}, \qquad q>{n\over r} \quad 
\hbox{if}\quad {n\over r}=p>1 .
$$ 
\medskip
\noindent Let $u\in W^{r,p}(\Omega)$. 
For any $x\in \Omega$ we set 
$$\Psi^x : y\in \Omega 
\to x+{{y-x}\over{\rho}},\qquad\Omega^*(x) =\Psi^x 
(\Omega(x))\,,$$ 

$$u^*=(u^x)^* : z\in \Omega^*(x) \longrightarrow u(x+\rho(z-x))\,.$$ 
We note that
$$
u^*\in W^{r,p}(\Omega^*(x)).
$$ 
We also note that, in consequence of $h_1)$ $\Omega^*(x)$ has 
the cone property, with the characteristic cone having height and 
opening independent of $x$. 

On the other hand, if $\tau={q/ p}$, from $h_2)$ we get 
$$\tau\geq 1,\qquad \tau>1 \quad \hbox{if}\quad {n\over r}={p>1},\qquad 
{{\tau-1}\over{p\tau}}\geq {{1\over p}-{r\over n}}.$$ 
From well-known imbedding theorems of Sobolev spaces (see, e.g., R.A. 
Adams [1]), we deduce that 
$$
u^*\in L^{{p\tau}\over {\tau-1}}(\Omega^*(x))
$$ 
and the following bound holds 
\begin{equation}\label{u*}
|u^*|_{{{p\tau}\over{\tau-1}},\Omega^*(x)}\leq c_0 
\|u^*\|_{W^{r,p}(\Omega^*(x))}, 
\end{equation} 
where $c_0=c_0(p,q,r,n)$ is a constant independent of $x$ and $u^*$. 

From (\ref{u*}) easily it follows that 
\begin{equation}\label{bound 1}
\rho^{-n\frac{(\tau-1)}{p\tau}}|u|_{{{p\tau}\over{\tau-1}},\Omega(x)} 
\le c_0
\sum_{|\alpha|\le r} \rho^{|\alpha|-{n\over p}}
|\partial^\alpha u|_{p,\Omega(x)}.
\end{equation} 

\bigskip\bigskip

\begin{thm}
If  $h_1)$ and $h_2)$  hold, 
then for any $g\in S^{q,\frac{s}{p}}(\Omega)$, $s\le p$, and for any $u\in W^{r,p}(\Omega)$ 
we get $gu\in L^p(\Omega)$ and 
\begin{equation}\label{bound 2} 
\|gu\|_{L^p(\Omega)}\le c\|g\|_{S^{q,\frac{s}{p}}(\Omega)} 
\|u\|_{W^{r,p}(\Omega)} ,
\end{equation}
where the constant $c=c(p,q,r,n)$ is independent of $g$ and $u$. 
\end{thm}

\medskip 

\begin{proof}
Let $u\in W^{r,p}(\Omega)$ and $g\in S^{q,\frac{s}{p}}(\Omega)$. 
By (\ref{map}) and by H{\" o}lder inequality it follows that 

\begin{equation}\label{bound 3}
\begin{split}
\int_\Omega |gu|^p dx &\le c_1 
\int_\Omega  \rho^{-n}  
\int_{\Omega(x)} |gu|^p dy\,dx\le 
c_1 \int_\Omega \rho^{-n}
\|g\|^p_{L^{p\tau}(\Omega(x))} 
\|u\|^p_{L^{\frac{p\tau}{\tau-1}}(\Omega(x))}dx \le
\\&
\le c_1 
\|g\|^p_{S^{q,\frac{s}{p}}(\Omega)}
\rho^{-\left(s+n\frac{\tau-1}{\tau}\right)}  
\int_\Omega 
\|u\|^p_{L^{\frac{p\tau}{\tau-1}}(\Omega(x))}dx.
\end{split}
\end{equation}
On the other hand from (\ref{bound 1}) and Lemma 4.2 we obtain 

\begin{equation}\label{bound 4}
\begin{split}
\rho^{-n\frac{(\tau-1)}{\tau}}\int_\Omega  
\|u\|^p_{L^\frac{p\tau}{\tau-1}(\Omega(x))} dx &\le 
 c_0 \, \sum_{|\alpha|\le r}\rho^{|\alpha|p-n}
\int_\Omega 
\|\partial^\alpha u\|^p_{L^p(\Omega(x))} dx \le
\\ &
\le c_1 \sum_{|\alpha|\le r} \rho^{|\alpha|p}
\|\partial^\alpha u\|^p_{L^p}.
\end{split}
\end{equation}
From (\ref{bound 3}), (\ref{bound 4}) the inequality (\ref{bound 2}) 
follows.

\end{proof}

\bigskip
The following theorem is a consequence of the embedding result stated in the\break
Theorem 4.3 (see result of Fefferman \cite{10} 
and also \cite{9} for a simplified proof).

\bigskip
\begin {thm}
If  $h_1)$ and $h_2)$  hold, 
then for any $g\in S^{q,\frac{s}{p}}(\Omega)$, $s\le p$,  
we get 
$$
\|gu\|_{L^p(\Omega)}\le c 
\|g\|_{S^{q,\frac{s}{p}}(\Omega)} 
\|\nabla u\|_{L^p(\Omega)} 
\qquad \forall u\in W^{1,p}_0(\Omega).
$$ 
where the constant $c=c(p,q,n)$ is independent of $g$ and $u$. 
\end{thm}

\medskip

\begin{proof}
Taking in mind the Hardy inequality,
the proof is a direct consequence of the Theorem 4.3 when $r=1$.
\end{proof}

\renewcommand\refname{{

\end{document}